\def\a{\alpha}
  \def\b{\beta}

  \def\e{\epsilon}

  \def\pf{{\it Proof. }$\;\;$}
  \def\hal{\unskip\nobreak\hfil\penalty50\hskip10pt\hbox{}\nobreak
  \hfill\vrule height 5pt width 6pt depth 1pt\par\vskip 2mm}

    \documentclass[11pt,a4paper]{article}
            \usepackage{amsfonts}

  \newtheorem{thm}{Theorem}[section]
  \newtheorem{prop}[thm]{Proposition}
  \newtheorem{lem}[thm]{Lemma}
  \newtheorem{cor}[thm]{Corollary}

\begin{document}

 \title{Simple groups, interleaved products, complexity and conjectures of Gowers and Viola}
  \author{Aner Shalev \\ Einstein Institute of Mathematics \\ The Hebrew University of Jerusalem
  \\ Jerusalem 91904 \\ Israel }

  \date{}
  \maketitle

\footnotetext{
The author acknowledges the support of an Israel Science Foundation grant 1117/13
and of the Vinik Chair of Mathematics which he holds. He also thanks the organizers
of the Conway Conference (November 2015) and Princeton University for their hospitality
while this work was carried out.}

\vspace{4mm}
\begin{abstract}
We study the distribution of products of conjugacy classes in finite simple groups, obtaining
various effective uniformity results, which give rise to an approximation to a conjecture
of Thompson.

Our results, combined with work of Gowers and Viola, also lead to the solution of recent conjectures
they posed on interleaved products and related complexity lower bounds, extending their work
on the groups SL$(2,q)$ to all (nonabelian) finite simple groups.

In particular it follows that, if $G$ is a finite simple group, and $A, B \subseteq G^t$
for $t \ge 2$ are subsets of fixed positive densities, then, as $a = (a_1, \ldots , a_t) \in A$
and $b = (b_1, \ldots , b_t) \in B$ are chosen uniformly, the interleaved product
$a \bullet b := a_1b_1 \cdots a_tb_t$ is almost uniform on $G$ (with quantitative estimates)
with respect to the $\ell_{\infty}$-norm.

It also follows that the communication complexity of an old decision problem related to
interleaved products of $a, b \in G^t$ is at least $\Omega(t \log |G|)$ when
$G$ is a finite simple group of Lie type of bounded rank, and at least
$\Omega(t \log \log |G|)$ when $G$ is any finite simple group. Both these bounds
are best possible.
\end{abstract}

\newpage

\section{Introduction}

The main purpose of this paper is to provide affirmative solutions to some
conjectures of Gowers and Viola -- see \cite{GV1, GV2, GV3}.
These papers contain interesting results in Group Theory (interleaved products)
and in Computer Science (complexity lower bounds) for the family of two-dimensional
special linear groups SL$(2,q)$.
Here we extend these results to all finite simple groups of Lie type of bounded Lie rank,
and in a weaker (yet quantitative) form to all finite simple groups.
In fact all our results here also apply (with similar proofs) to all finite
quasisimple groups, namely finite perfect groups $G$ such that $G/Z(G)$ is simple.

Throughout this paper simple groups are taken to be nonabelian, and we assume the
Classification of finite simple groups.
Since our results are of asymptotic nature we may ignore the sporadic groups and restrict
our attention to simple groups of Lie type and to alternating groups $A_n$.

Our main contribution is Theorem 1.1 below, on the distribution of products of elements from two
random conjugacy classes of a finite simple group. This quantitative result
has a number of direct consequences -- see results 1.2-1.5 below; see also \cite{Sh1, Sh2} for
earlier results in this direction, which are not sufficiently strong for the current applications.
In particular we derive (in Theorem 1.4 below) a quantitative approximation to a conjecture of Thompson
(see \cite{AH} and \cite{EG}) which is still open for simple groups of Lie type over tiny fields.
The combination of Corollary 1.5 with reductions and statements from \cite{GV3} yields
various applications to interleaved products and complexity, some of which are mentioned
briefly in Sections 1 and 3 of this paper.

We start with some notation which we will use throughout this paper.
Let $G$ be a finite group and let $x, y, g \in G$. Let
$p_{x,y}(g)$ denote the probability that $g = x' y'$, where $x'$ is
a random conjugate of $x$ and $y'$ is a random conjugate of $y$
(with respect to the uniform distribution).
Then $p_{x,y}$ is a probability distribution on $G$. Let $||p_{x,y}||_2^2$
denote the square of its $\ell_2$-norm, namely
\[
||p_{x,y}||_2^2 = \sum_{g \in G} p_{x,y}(g)^2.
\]
By $IrrG$ we denote the set of complex irreducible characters of $G$.
We define the Witten zeta function $\zeta_G$ of $G$ by
\[
\zeta_G(s) = \sum_{\chi \in IrrG} \chi(1)^{-s},
\]
where $s$ is a real number. This function plays a key role in our proofs.
Our notation for finite simple groups of Lie type, their rank and their
underlying field, follows that of \cite{LiSh3}.

Our main theorem below implies that for finite simple groups $G$,
and for almost all $x, y \in G$, the distribution $p_{x,y}$
is very close to uniform in the $\ell_2$ sense. For the applications
we prove a rather general quantitative result, where $x, y$ need not be independent.

\begin{thm} Let $G$ be a finite simple group.
Let $\nu$ be a probability distribution on $G^2$ which
projects to uniform distributions on each coordinate.
Choose $(x, y) \in G^2$ according to the distribution $\nu$
(so that $x$ is uniform in $G$ and so is $y$, but they are not
assumed to be independent).

(i) If $G = A_n$ then, for some absolute constant $c$, the $\nu$-probability that
$||p_{x,y}||_2^2 \le |G|^{-1}(1 + c n^{-2/3})$ is greater than $1 - c n^{-2/3}$.

(ii) For any $\e > 0$ there is $r(\e)$ such that if $r \ge r(\e)$ and
$G$ is a group of Lie type of rank $r$ over the field with $q$ elements,
then the $\nu$-probability that
$||p_{x,y}||_2^2 \le |G|^{-1}(1 + q^{-(2/3-\e)r})$ is greater than $1 - q^{-(2/3-\e)r}$.

(iii) If $G$ is a group of Lie type of rank $r$, then there exists $c = c(r) > 0$
such that the $\nu$-probability that
$||p_{x,y}||_2^2 \le |G|^{-1}(1 + |G|^{-c})$ is at least $1 - |G|^{-c}$.

\end{thm}

We can also show that if $G$ is alternating or a group of
Lie type of unbounded rank then part (iii) above does not
hold for an absolute constant $c > 0$.

Theorem 1.1 is a particular case of more general results, which also yield better
bounds on $||p_{x,y}||_2^2$ (possibly with lower probabilities) -- see Theorems
2.4 and 2.6 below. In particular we show that if $G$ is a finite simple group of Lie type of rank $r$
over the field with $q$ elements, then the probability that
\[
||p_{x, y}||_2^2 \le  1 + q^{-(2-\e)r}
\]
is at least $1 - q^{-\frac{1}{3}\e r}$, for any $\e > 0$ and $r \ge r(\e)$ (see Corollary 2.7 for this
and for a similar result for alternating groups).

Theorem 1.1 and its variants have several interesting consequences which we discuss below.

Let $U$ be the uniform distribution on $G$. A trivial calculation shows that
the $\ell_2$-distance between the distributions $p_{x,y}$ and $U$ satisfies
\[
||p_{x,y}-U||_2^2 = ||p_{x,y}||_2^2 - |G|^{-1},
\]
which can be effectively bounded (for almost all $x, y$) by Theorem 1.1 above.

Next, consider the $\ell_1$-distance (also known as the statistical distance, or the total variation distance
up to normalization) between the distributions $p_{x,y}$ and $U$, defined by
\[
||p_{x,y} - U||_1 = \sum_{g \in G} |p_{x,y}(g) - |G|^{-1}|.
\]
In \cite[2.5]{Sh1} it is shown that if $G$ is a finite simple group, and $x, y \in G$
are chosen uniformly and independently, then we have $||p_{x,y} - U||_1 = o(1)$ with probability
at least $1 - o(1)$, where, throughout this paper, $o(1)$ is a quantity tending to $0$ as
$|G| \rightarrow \infty$. Here we obtain a stronger result, where $x,y$ need not be independent, and the
estimates are effective and close to best possible.

\begin{cor} Let $G$ be a finite simple group.
Let $\nu$ be a probability distribution on $G^2$ which
projects to uniform distributions on each coordinate.
Choose $(x, y) \in G^2$ according to the distribution $\nu$

(i) If $G = A_n$ then, for some absolute constant $c$, the $\nu$-probability that
$||p_{x,y}-U||_1 \le c n^{-1/3}$ is greater than $1 - c n^{-2/3}$.

(ii) For any $\e > 0$ there is $r(\e)$ such that if $r \ge r(\e)$ and
$G$ is a group of Lie type of rank $r$ over the field with $q$ elements,
then the $\nu$-probability that
$||p_{x,y}-U||_1 \le q^{-(1/3-\e)r})$ is greater than $1 - q^{-(2/3-\e)r}$.

(iii) If $G$ is a group of Lie type of rank $r$, then there exists $c = c(r) > 0$
such that the $\nu$-probability that
$||p_{x,y}-U||_1 \le |G|^{-c}$ is at least $1 - |G|^{-2c}$.
\end{cor}

Corollary 1.2 follows easily from Theorem 1.1. Indeed, by the Cauchy-Schwarz inequality
we have
\[
||p_{x,y}-U||_1 \le ||p_{x,y} - U||_2 \cdot |G|^{1/2} =  (||p_{x,y}||_2^2 - |G|^{-1})^{1/2}|G|^{1/2}.
\]
This means that, if $||p_{x,y}||_2^2 \le |G|^{-1}(1 + \delta)$ (where $\delta$ is given
by Theorem 1.1), then $||p_{x,y}-U||_1 \le \delta^{1/2}$.

Part (iii) of Corollary 1.2 above extends \cite[1.10]{GV2} and \cite[1.12]{GV3} dealing with
$G = $SL$(2,q)$ to all finite simple groups of Lie type of bounded rank.

Another application of Theorem 1.1 concerns the size of the product $x^Gy^G$ of the
conjugacy classes of $x$ and of $y$ in $G$. A famous conjecture of J.G. Thompson states
that every finite simple group has a conjugacy class $x^G$ satisfying $(x^G)^2 = G$.
This was confirmed for alternating groups $A_n$ and for groups of Lie type over fields
with more than $8$ elements, see \cite{EG} and the references therein. However, the case of
classical groups over tiny fields remains open. See also \cite[1.1.4]{LST} and \cite[1.4]{GM}
for variations on Thompson's conjecture, dealing with products $x^G y^G$ of two conjugacy classes.

The following quantitative result shows that $x^Gy^G$ usually almost covers $G$;
this applies also to $(x^G)^2$, since $x,y$ need not be independent, so we may take $x=y$.

\begin{cor} Let $G$ be a finite simple group.
Let $\nu$ be a probability distribution on $G^2$ which
projects to uniform distributions on each coordinate.
Choose $(x, y) \in G^2$ according to the distribution $\nu$

(i) If $G = A_n$ then, for some absolute constant $c$, the $\nu$-probability that
$|x^Gy^G| \ge (1- c n^{-2/3})|G|$ is greater than $1 - c n^{-2/3}$.

(ii) For any $\e > 0$ there is $r(\e)$ such that if $r \ge r(\e)$ and
$G$ is a group of Lie type of rank $r$ over the field with $q$ elements,
then the $\nu$-probability that
$|x^Gy^G| \ge (1-q^{-(2/3-\e)r})|G|$ is greater than $1 - q^{-(2/3-\e)r}$.

(iii) If $G$ is a group of Lie type of rank $r$, then there exists $c = c(r) > 0$
such that the $\nu$-probability that
$|x^Gy^G| \ge (1-|G|^{-c})|G|$ is at least $1 - |G|^{-c}$.
\end{cor}

To prove this, note that if $|x^Gy^G| = (1-\delta)|G|$ then
$||p_{x,y}||_2^2 \ge (1-\delta)^{-1}|G|^{-1} \ge (1+ \delta)|G|^{-1}$.
Corollary 1.3 now follows immediately from Theorem 1.1.

By a remark following Theorem 1.1 it also follows that,
if $G$ is a finite simple group of Lie type of rank $r$ over the field with $q$ elements,
then the probability that $|x^G y^G| \ge  (1 - q^{-(2-\e)r})|G|$ is at least $1 - q^{-\frac{1}{3}\e r}$,
for any $\e > 0$ and $r \ge r(\e)$.
This gives rise to the following approximation to Thompson's conjecture in the open case of
classical groups over tiny fields.

\begin{thm}
For any $\e > 0$ there is $r(\e)$ such that if $r \ge r(\e)$ and $G$ is a finite simple group
of Lie type of rank $r$ over the field with $q$ elements, then there exists a conjugacy class
$x^G$ of $G$ such that $|(x^G)^2| \ge (1-q^{-(2-\e)r})|G|$.
\end{thm}

Theorem 1.1 applies in various situations; these include the cases where $x, y$ are uniform and independent,
when $x$ is uniform and $y=x$, and more generally, when $x$ is uniform and $y = f(x)$,
where $f:G \rightarrow G$ is any fixed bijection.

In particular, if we fix $a \in G$ and let $f$ be the bijection sending $x$ to $x^{-1}a$,
we obtain the following.

\begin{cor} Let $G$ be a finite simple group, let $a \in G$ be any fixed element,
let $x \in G$ distribute uniformly over $G$ and let $y = x^{-1}a$. Then
$p_{x,y}$ satisfies the conclusions (i)-(iii) of Theorem 1.1.
\end{cor}

In the case of $G = $ SL$(2,q)$ this result is proved in \cite[1.13]{GV3}.
It is also stated in \cite{GV3} that if Corollary 1.5 above holds for a family
of finite groups $G$ then these groups satisfy a variety of interesting results,
proven earlier only for SL$(2,q)$. We mention now briefly some of these applications,
while some more will be discussed in Section 3.

Recall that for a group $G$, a positive integer $t \ge 2$, and two $t$-tuples
$a = (a_1, \ldots , a_t), b = (b_1, \ldots , b_t) \in G^t$, the {\em interleaved product}
$a \bullet b$ of $a$ and $b$ is defined by
\[
a \bullet b = a_1b_1a_2b_2 \cdots a_tb_t \in G.
\]
The density of a subset $A \subseteq G^t$ is defined by $|A|/|G|^t$.

\newpage

\begin{thm} Let $G$ be a finite simple group and $t \ge 2$ an integer.
Let $A, B \subseteq G^t$ be subsets of positive densities $\a$ and
$\b$ respsectively.
If $a$ and $b$ are selected uniformly from $A$ and $B$, then, for each
$g \in G$, the probability that $a \bullet b = g$ is of the form
$(1 + o(1))|G|^{-1}$.

In particular, if $G$ is sufficiently large (given $\a$ and $\b$),
then $A \bullet B = G$.
\end{thm}

Thus $a \bullet b$ (for $a \in A$ and $b \in B$) is almost uniformly distributed
in the $\ell_{\infty}$-norm.

Theorem 1.6 above follows from stronger bounds as follows.
Let $\alpha = |A|/|G|^t$ and $\beta = |B|/|G|^t$ be the densities of
$A$ and $B$ respectively.
If the simple group $G$ above is of Lie type of bounded rank then
we obtain
\[
|Prob(a \bullet b = g) - |G|^{-1}| \le (\alpha \beta)^{-1} |G|^{-1 - ct},
\]
where $c > 0$ depends only on the rank of $G$. This extends Theorem 1.7 of \cite{GV2}
(which is Theorem 1.8 of \cite{GV3}) dealing with SL$(2,q)$.

If $G$ is any simple group of Lie type of rank $r$ (which is not necessarily bounded)
we obtain
\[
|Prob(a \bullet b = g) - |G|^{-1}|  \le (\alpha \beta)^{-1} q^{-c r t} |G|^{-1},
\]
where $c > 0$ is an absolute constant.

Finally, if $G = A_n$ then, for some absolute positive constant $c$ we have
\[
|Prob(a \bullet b = g) - |G|^{-1}|  \le (\alpha \beta)^{-1} n^{-c t} |G|^{-1}.
\]

These results generalize the case when the subsets $A, B$ are product sets, and the
related distribution can then be analyzed using Gowers' paper \cite{G} and the paper
\cite{BNP} by Babai, Nikolov and Pyber.

Applications of Corollary 1.5 to certain complexity lower bounds and related
conjectures of Gowers and Viola will be described in Section 3 below.
In fact Corollary 1.5 also extends additional results from \cite{GV1, GV2, GV3},
and is likely to have further applications in subsequent works.

We note that while the proofs in \cite{GV1, GV2, GV3} avoid representation theory, we use
it as our main tool, which sometimes yields shorter proofs of more general results.

\medskip

{\bf Acknowledgment.}
I am grateful to Tim Gowers for interesting conversations, for sending me the preprint \cite{GV3}
and for asking me about possible extensions to other simple groups.

\newpage

\section{Proof of Theorem 1.1}

This section is devoted to the proof of Theorem 1.1.
In fact we prove somewhat stronger results (see Theorems 2.4 and 2.6 below)
from which Theorem 1.1 follows.
We need some preparations.

\begin{lem} Let $G$ be a finite group, and $x, y \in G$.
Then we have
\[
||p_{x,y}||_2^2 = |G|^{-1} \sum_{\chi \in IrrG}|\chi(x)|^2|\chi(y)|^2/\chi(1)^2.
\]
\end{lem}

\pf It is well known that
\[
p_{x,y}(g) = |G|^{-1} \sum_{\chi \in IrrG} \chi(x)\chi(y)\chi(g^{-1})/\chi(1).
\]
Therefore
\[
||p_{x,y}||_2^2 = |G|^{-2}\sum_{g \in G} [\sum_{\chi \in IrrG} \chi(x)\chi(y)\chi(g^{-1})/\chi(1)]^2.
\]
This yields
\[
||p_{x,y}||_2^2 = |G|^{-2} \sum_{g\in G} \sum_{\chi, \psi \in IrrG} \chi(x)\chi(y)\psi(x)\psi(y)/(\chi(1)\psi(1))
\cdot \chi(g^{-1}) \psi(g^{-1}).
\]
Changing the order of summation we obtain
\[
||p_{x,y}||_2^2 = |G|^{-2} \sum_{\chi, \psi \in IrrG} \chi(x)\chi(y)\psi(x)\psi(y)/(\chi(1)\psi(1))
\cdot \sum_{g \in G} \chi(g^{-1}) \psi(g^{-1}),
\]
which, by the orthogonality relations, vanishes unless $\psi = \overline{\chi}$, yielding
\[
||p_{x,y}||_2^2 = |G|^{-1} \sum_{\chi \in IrrG}|\chi(x)|^2|\chi(y)|^2/\chi(1)^2.
\]
\hal
\medskip

\begin{prop} Let $G$ be a finite simple group. Then

(i) For a fixed real number $s > 1$ we have $\zeta_G(s) = 1 + o(1)$.

(ii) If $G$ is a group of Lie type and $s>1$ then
there exists $c>0$ depending only on $s$ and on the rank of $G$
such that $\zeta_G(s) \le 1 + |G|^{-c}$.

(iii) If $G = A_n$ then for any fixed real number $s > 0$ we have
$\zeta_G(s) = 1 + O(n^{-s})$.

(iv) For any fixed real numbers $s, \e > 0$ there is a number
$r(s,\e)$ such that, if $G$ is a group of Lie type of rank $r \ge r(s,\e)$ over
a field with $q$ elements, then we have
$\zeta_G(s) \le 1 + q^{-(s-\e)r}$.

\end{prop}

\pf Parts (i) and (iii) are proved in \cite[2.7]{LiSh1} for alternating groups.

The proof of part (i) for groups of Lie type (and in fact for all finite
quasisimple groups) is given in \cite[1.1]{LiSh2}.

To prove part(ii), let $G$ be a group of Lie type of rank $r$
over the field with $q$ elements. Let $k(G)$ be the number
of conjugacy classes of $G$. It is known (see \cite[1.1]{FG}) that
$k(G) \le c_1 q^r$ for some absolute constant $c_1$.
It is also known (see \cite{LS}) that there is an absolute
constant $c_2 > 0$ such that $\chi(1) \ge c_2 q^r$ for every
nontrivial character $\chi \in IrrG$.
It follows that, for $s>1$,
\[
\zeta_G(s) \le 1 + c_1 q^r (c_2 q^r)^{-s} \le 1  + c_3 q^{-r(s-1)},
\]
where $c_3$ depends on $s$. Since $|G| \le q^{4r^2}$ this yields
\[
\zeta_G(s) \le 1+ |G|^{-c},
\]
where $c$ depends on $s$ and $r$.

The proof of part (iv) applies a variation on arguments from \cite{LiSh3}.
First note that it suffices to prove part (iv) for classical groups
of large rank (since we may choose $r(s,\e)$ large enough).
In the proof of Theorem 1.2 of \cite{LiSh3} it is shown
that, for every $c > 0$ there exist $N=N(c)$ and $c_2 = c_2(c)$ such that if the classical
group $G$ has a natural module of dimension $n \ge N$, and $s > 0$, then
\[
\zeta_G(s) \le 1 + c_1 q^c c_2^{\sqrt{n}} q^{-s(n-1)/2} + c_3 c_4^{-s} q^n q^{-csn},
\]
where $c_1, c_3, c_4$ are absolute constants. Examination of the arguments there
shows that the term $q^{-s(n-1)/2}$ may be replaced by $q^{-sr}$ where
$r$ is the rank of $G$.

Now, given $s$, set $c = 2/s + 1$ (rather than $c = 2/s$ as in \cite{LiSh3}).
Then, for $n \ge N(c)$ we have
\[
\zeta_G(s) \le 1 + c_1 q^c c_2^{\sqrt{n}} q^{-sr} + c_3 c_4^{-s} q^{-(s+1)n}.
\]
Since $r \le n$ it easily follows (focusing on the dominant terms) that for any
$\e>0$ there exists $r(s, \e) \ge N(c)$ such that for $r \ge r(s, \e)$ we have
\[
\zeta_G(s) \le 1 + q^{-(s-\e)r}.
\]
This completes the proof.

\hal
\medskip

We note that parts (iii) and (iv) above are almost best possible, since
they show that $\zeta_G(s)$ is well approximated by its two first summands.

\begin{prop} Let $G$ be a finite simple group of Lie type of rank $r$
over a field with $q$ elements.

(i) There is a constant $c > 0$ depending only on $r$,
such that, if $x, y$ distribute uniformly over $G$ (but may be
dependent), then
\[
\sum_{\chi \in IrrG}|\chi(x)|^2|\chi(y)|^2/\chi(1)^2 \le 1 + |G|^{-c}
\]
holds with probability at least $1 - |G|^{-c}$.

(ii) There is an absolute constant $c > 0$ and a constant $c'$ depending on $r$,
such that, if $G \not\in S$, where
\[
S = \{ L_2(q), L_3^{\pm}(q), L_4^{\pm}(q), D_4^{\pm}(q), D_5^{\pm}(q) \},
\]
and $x, y$ distribute uniformly over $G$ (but may be dependent), then
\[
\sum_{\chi \in IrrG}|\chi(x)|^2|\chi(y)|^2/\chi(1)^2 \le 1 + c' q^{-(2r-1)}
\]
holds with probability at least $1 - c/q$.

\end{prop}

\pf Let $G$ be of rank $r$ over the field with $q$ elements.
It is known that the probability that $x \in G$ is regular semisimple is at
least $1-c_1/q$ for an absolute constant $c_1 > 0$ (see \cite{GL} for a more detailed result).
Therefore the probability that both $x$ and $y$ are regular semisimple is at least $1-2c_1/q$.
Note that the this also holds if $x, y$ are dependent.

If $x, y \in G$ are regular semisimple then $|\chi(x)|, |\chi(y)| \le b$,
where $b$ depends only on the rank of $G$ (see e.g. \cite[4.4]{Sh2}). This yields
\[
\sum_{1 \ne \chi \in IrrG}|\chi(x)|^2|\chi(y)|^2/\chi(1)^2  \le b^4 (\zeta(2) - 1).
\]
Part (i) now follows using Proposition 2.2(ii).

To prove part (ii) we shall derive better bounds than those provided by the
proof of Proposition 2.2(ii). Note that this proof yields $\zeta_G(2) \le 1 + c q^{-r}$
for some absolute constant $c$. Our main tool is Proposition 6.2 of \cite{LiSh3}.
It shows that if $G$ is a finite simple group of Lie type over a field with $q$ elements,
and $G \not\in S$, and $\chi \in Irr G$ is not a Weil character, then $\chi(1) > \max(q^{3r/2}, q^{2r-3})$.
Now, the number of Weil characters is at most $\max(q+1, 4) \le q+2$, and this yields (for $G \not\in S$)
\[
\zeta_G(2) \le 1 + (q+2)(c_2q^r)^{-2} + c_1q^r(q^{3r/2})^{-2},
\]
where $c_1, c_2$ are as in the proof of Proposition 2.2(ii). It now easily follows
that, for some absolute constant $c_4 > 0$ we have
\[
\zeta_G(2) \le 1 + c_4 q^{-(2r-1)}.
\]
This, combined with the proof of part (i), yields part (ii) of the proposition.

\hal
\medskip

We can now derive the main result of this section for groups of bounded rank.

\begin{thm} Let $G$ be a finite simple group of Lie type of rank $r$
over a field with $q$ elements.

(i) There is a constant $c > 0$ depending only on $r$,
such that, if $x, y$ distribute uniformly over $G$ (but may be
dependent), then
\[
||p_{x,y}||_2^2 \le 1 + |G|^{-c}
\]
holds with probability at least $1 - |G|^{-c}$.

(ii) There is an absolute constant $c > 0$ and a constant $c'$ depending only on $r$,
such that, if $G \not\in S$, where
\[
S = \{ L_2(q), L_3^{\pm}(q), L_4^{\pm}(q), D_4^{\pm}(q), D_5^{\pm}(q) \},
\]
and $x, y$ distribute uniformly over $G$ (but may be dependent), then
\[
||p_{x,y}||_2^2 \le 1 + c' q^{-(2r-1)}
\]
holds with probability at least $1 - c/q$.

\end{thm}

\pf This follows from Lemma 2.1 and Proposition 2.3 above.

\hal
\medskip

We now turn to alternating groups and groups of Lie type of unbounded rank.

\begin{prop} Let $G$ be a finite simple group. Let $x, y$ distribute
uniformly over $G$ (but they may be dependent). Fix $s$ with $s > 0$.

(i) If $G = A_n$ then for some absolute constant $c$ the probability that
\[
\sum_{\chi \in IrrG}|\chi(x)|^2|\chi(y)|^2/\chi(1)^2 \le 1 + c n^{-(2-2s)}
\]
is at least $1- c n^{-s}$.

(ii) If $G$ is a finite simple group of Lie type of rank $r$ over the field with $q$
elements, then the probability that
\[
\sum_{\chi \in IrrG}|\chi(x)|^2|\chi(y)|^2/\chi(1)^2 \le  1 + q^{-(2 - 2s -\e)r}
\]
is at least $1 - q^{-(s-\e)r}$, for any $\e > 0$ and $r \ge r(s,\e)$.
\end{prop}

\pf
It follows from \cite[2.2]{Sh1} that, for any finite group $G$, a fixed $s > 0$ and a uniformly
distributed $x \in G$, the probability that
\[
|\chi(x)| \le \chi(1)^{s/2}
\]
for all $\chi \in IrrG$ is greater than $2 - \zeta_G(s) = 1 - (\zeta_G(s)-1)$.

We conclude that for uniform (possibly dependent) $x, y \in G$,
the probability that $|\chi(x)| \le \chi(1)^{s/2}$ and $|\chi(y)| \le \chi(1)^{s/2}$
for all $\chi \in IrrG$ is greater than $3 - 2 \zeta_G(s) = 1 - 2(\zeta_G(s)-1)$.
Hence the inequality
\[
\sum_{\chi \in IrrG}|\chi(x)|^2|\chi(y)|^2/\chi(1)^2 \le \zeta_G(2-2s)
\]
holds with probability greater than $3-2 \zeta_G(s)$.

We now apply Proposition 2.2. If $G = A_n$ then by part (iii) of this result we have
\[
\zeta_G(2-2s) \le 1 + c n^{-(2-2s)} \;\; {\rm and} \; \; 3-2 \zeta_G(s) \ge 1 - 2c n^{-s}
\]
where $c$ is an absolute constant. Plugging this in the previous probability estimate
(replacing $c$ by $c/2$) proves part (i).

Now let $G$ be a group of Lie type of rank $r$ over the field with $q$ elements. Then
part (iv) of Proposition 2.2 yields
\[
\zeta_G(2-2s) \le 1 + q^{-(2-2s-\e)r} \;\; {\rm and} \; \; 3-2 \zeta_G(s) \ge 1 - 2 q^{-(s-\e)r}
\]
for any $\e >0$ and $r \ge r(s,\e)$.

Part (ii) follows from this and the above discussion (replacing $\e$, say, by $\e/2$).

\hal
\medskip

We can now prove the main result of this section for groups of unbounded rank.

\begin{thm} Let $G$ be a finite simple group. Let $x, y$ distribute
uniformly over $G$ (but they may be dependent). Fix $s$ with $s > 0$.

(i) If $G = A_n$ then for some absolute constant $c$ the probability that
\[
||p_{x,y}||_2^2 \le 1 + c n^{-(2-2s)}
\]
is at least $1- cn^{-s}$.

(ii) If $G$ is a finite simple group of Lie type of rank $r$ over the field with $q$
elements, then the probability that
\[
||p_{x,y}||_2^2 \le  1 + q^{-(2 - 2s -\e)r}
\]
is at least $1 - q^{-(s-\e)r}$, for any $\e > 0$ and $r \ge r(s,\e)$.
\end{thm}

\pf

This follows immediately from Lemma 2.1 and Proposition 2.5 above.

\hal
\medskip

{\bf Proof of Theorem 1.1:} Parts (i) and (ii) of the theorem follow from
Theorem 2.6 above by substituting $s = 2/3$.
Part (iii) of the theorem is part (i) of Theorem 2.4 above.

\medskip

The following consequence of Theorem 2.6 will also be useful.

\begin{cor} Let $G$ be a finite simple group. Let $x, y$ distribute
uniformly over $G$ (but they may be dependent).

(i) If $G = A_n$ then for any $\e > 0$ there exists $n(\e)$ such that
for any $n \ge n(\e)$ the probability that
\[
||p_{x, y}||_2^2 \le 1 + n^{-(2-\e)}
\]
is at least $1- n^{-\e/3}$.

(ii) If $G$ is a finite simple group of Lie type of rank $r$ over the field with $q$
elements, then the probability that
\[
||p_{x, y}||_2^2 \le  1 + q^{-(2-\e)r}
\]
is at least $1 - q^{-\frac{1}{3}\e r}$, for any $\e > 0$ and $r \ge r(\e)$.
\end{cor}

\pf

Part (i) follows from part (i) of  Theorem 2.6 with $s = 2\e/5$.

To prove part (ii) apply part (ii) of Theorem 2.6 with $s = 4 \e$.
We see that $||p_{x,y}||_2^2 \le 1 + q^{(2-9 \e)r}$ with probability
at least $1 - q^{-3 \e r}$. Replacing $\e$ by $\e/9$ we obtain the
result.

\hal
\medskip

Note that, if we replace $3$ in the conclusions of Corollary 2.7 by
any fixed number greater than $2$, the conclusions will still hold.

\section{Complexity applications}

In this section we briefly describe applications of our main results
to complexity lower bounds related to interleaved products. We follow
definitions and statements from \cite{GV1, GV2, GV3}.

Consider the following promise problem introduced in 1984 in \cite{ESY}.
Let $G$ be a finite group and $t \ge 2$ an integer. Suppose Alice receives
a $t$-tuple $a \in G^t$ and Bob receives a $t$-tuple $b \in G^t$.
Suppose we are promised that the interleaved product $a \bullet b \in G$ is
one of two given elements $g, h \in G$. The task of Alice and Bob is to
decide whether $a \bullet b = g$ or $a \bullet b = h$. What can we say about
the communication complexity of this problem?

Recall that $O(n)$ denotes numbers bounded above by $cn$ for some
constant $c$, while $\Omega(n)$ denotes numbers bounded below
by $cn$ for some positive constant $c$.

Note that a trivial upper bound for the communication complexity above is $O(t \log |G|)$.
It is shown in \cite{GV1, GV2, GV3} that this upper bound is tight for
$G = $SL$(2,q)$, namely, in this case the communication complexity
is at least $\Omega(t \log |G|)$. Corollary 1.5 combined with
reductions and statements from \cite{GV3} extend this as follows.

\begin{thm}
The above communication complexity is at least $\Omega(t \log |G|)$
whenever $G$ is a finite simple group of Lie type of bounded rank.
\end{thm}

For general finite simple groups we obtain the following.

\begin{thm} The above communication complexity is at least
 $\Omega(t \log \log |G|)$ whenever $G$ is a finite simple group.
 If $G$ is a finite simple group of Lie type, then the communication
 complexity is at least $\Omega(t \sqrt{\log |G|})$.
 \end{thm}

 The first assertion in Theorem 3.2 was conjectured by Gowers and Viola
 (see \cite{GV1, GV2, GV3}). This complexity lower bound is tight for
 alternating groups (see \cite{MV}).

The next result easily implies the complexity bounds in Theorems 3.1 and 3.2;
it extends Theorem 1.2 of \cite{GV2, GV3} which deals with $G=$ SL$(2,q)$.

\begin{thm} Let $G$ be a finite simple group and let $t \ge 2$ be an integer.
Let $P: G^t \times G^t \rightarrow \{ 0, 1 \}$ be a (randomized public-coin)
$c$-bit communication protocol. For $g \in G$ let $p_g$ denote the probability
that $P(a,b)=1$ assuming $a \bullet b = g$. Then for any $g, h \in G$ we have

(i) $|p_g - p_h| \le 2^c |G|^{-\Omega(t)}$ if $G$ is a group of Lie type
of bounded rank.

(ii) $|p_g - p_h| \le 2^c q^{-\Omega(rt)}$ if $G$ is a group of Lie type
of rank $r$.

(iii) $|p_g - p_h| \le 2^c n^{-\Omega(t)}$ if $G = A_n$.
\end{thm}

This result follows from Corollary 1.5 combined with statements from \cite{GV3}.

The following is an immediate consequence of Theorem 3.3.

\begin{cor} With the above notation we have $|p_g - p_h| \le 2^c (\log |G|)^{-\Omega(t)}$
for all finite simple groups $G$.
\end{cor}

This proves Conjecture 1.3 in \cite{GV2, GV3}.

\end{document}